\documentclass[a4paper,11pt]{article}
\usepackage{amsmath,amssymb,amsthm,fullpage,tikz,authblk}

\newtheorem{theorem}{Theorem}[section]
\newtheorem{corollary}[theorem]{Corollary}
\newtheorem{conjecture}[theorem]{Conjecture}
\newtheorem{question}[theorem]{Question}
\newtheorem{proposition}[theorem]{Proposition}

\theoremstyle{remark}
\newtheorem*{remark}{Remark}

\newcommand{\z}{\mathbb{Z}}
\newcommand{\zp}{\mathbb{Z}_+}
\newcommand*{\abs}[1]{\lvert #1\rvert}
\newcommand*{\babs}[1]{\bigl| #1\bigr|}
\newcommand*{\ceil}[1]{\lceil #1\rceil}
\newcommand*{\bceil}[1]{\bigl\lceil #1\bigr\rceil}

\newcommand*{\sm}[1][G]{\operatorname{sum}(#1)}
\newcommand*{\df}[1][G]{\operatorname{dif}(#1)}
\newcommand*{\smb}[1][\hat{K}_n^{(k,k)}]{\operatorname{sum}\bigl(#1\bigr)}
\newcommand*{\dfb}[1][\hat{K}_n^{(k,k)}]{\operatorname{dif}\bigl(#1\bigr)}

\title{Sum index, difference index and exclusive sum number of graphs}
\author{John Haslegrave\thanks{Research supported by the UK Research and Innovation Future Leaders Fellowship MR/S016325/1.}}
\affil{Mathematical Institute, University of Oxford}

\begin{document}
\maketitle

\begin{abstract}We consider two recent conjectures of Harrington, Henninger-Voss, Karhadkar, Robinson and Wong concerning relationships between the sum index, difference index and exclusive sum number of graphs. One conjecture posits an exact relationship between the first two invariants; we show that in fact the predicted value may be arbitrarily far from the truth in either direction. In the process we establish some new bounds on both the sum and difference index. The other conjecture, that the exclusive sum number can exceed the sum index by an arbitrarily large amount, follows from known, but non-constructive, results; we give an explicit construction demonstrating it.

Simultaneously with the first version of this paper appearing, Harrington et al.\ updated their preprint with two counterexamples to the first conjecture; however, their counter\-examples only give a discrepancy of $1$, and only in one direction. They therefore modified the conjecture from an equality to an inequality; our results show that this is still false in general.
\end{abstract}

\section{Introduction}\label{intro}
Let $G$ be a finite simple graph and let $f:V(G)\to\z$ be an injective labelling. We may derive two natural edge labellings $f_+:E(G)\to \z$ and $f_-:E(G)\to\z$ from $f$ as follows.
\begin{align*}
f_+:uv&\mapsto f(u)+f(v);\\
f_-:uv&\mapsto \abs{f(u)-f(v)}.
\end{align*}
We consider two parameters of $G$ based on these labellings. The \textit{sum index}, $\sm$, introduced by Harrington and Wong \cite{HW20}, is the minimum number of different values taken by $f_+$ for any choice of $f$, and the \textit{difference index}, $\df$, is defined likewise with respect to $f_-$. Note that restricting $f$ to the positive integers would make no difference, since adding a constant to $f$ does not change the number of values taken by $f_+$ or $f_-$. In a recent preprint \cite{HHKRW}, Harrington, Henninger-Voss, Karhadkar, Robinson and Wong introduced the difference index and made the following bold conjecture, based on computations of exact values for several families of graphs.
\begin{conjecture}[{\cite[Conjecture 4.2]{HHKRW}}]\label{conj:sum-dif}Every nonempty simple graph $G$ satisfies \[\df=\ceil{\sm/2}.\] 
\end{conjecture}
They further supported Conjecture \ref{conj:sum-dif} by proving a one-directional version for bipartite graphs: if $G$ is bipartite then $\df\geq\ceil{\sm/2}$.

We show in Sections \ref{high-sum} and \ref{low-sum} that Conjecture \ref{conj:sum-dif} is false in general, and in fact the conjectured value can be arbitrarily far from the truth in either direction. In one direction, we prove this in the strongest possible form: we show that the sum index cannot be bounded by any function of the difference index, and in fact even for difference index $2$ there are connected graphs with arbitrary sum index. In the other direction, while the discrepancy from the conjectured value may be arbitrarily large, it is small relative to the size of the two indices. We give some new bounds on the sum and difference indices, and as a consequence we obtain the exact sum index of the prism graphs, which was not previously known.

In an updated version \cite{HHKRW2}, which appeared as a preprint simultaneously with the first version of this paper, Harrington et al.\ give two counterexamples to Conjecture \ref{conj:sum-dif2} in the second direction, although both of these only achieve a discrepancy of $1$. They therefore modify Conjecture \ref{conj:sum-dif} as follows.
\begin{conjecture}[{\cite[Conjecture 4.4]{HHKRW2}}]\label{conj:sum-dif2}Every nonempty simple graph $G$ satisfies \[\ceil{\sm/2}\leq\df\leq\sm.\]\end{conjecture}
Our results in Section \ref{high-sum} show that Conjecture \ref{conj:sum-dif2} is still false in a strong sense.

We also consider the exclusive sum number of a graph, which is closely related to the well-studied sum number, introduced by Harary \cite{Har90}. Harary defined a graph $G$ to be a \textit{sum graph} if there is an injective function $f:V(G)\to\zp$ such that for any distinct vertices $u,v\in V(G)$ we have $uv\in E(G)$ if and only if $f(u)+f(v)\in f(V(G))$; this $f$ is a \textit{sum labelling} of $G$. Necessarily any sum graph has at least one isolated vertex, since the vertex $v$ maximising $f(v)$ cannot have any neighbours; the \textit{sum number}, $\sigma(G)$, of a graph $G$ is the smallest number of isolated vertices which can be added to create a sum graph. Harary asked for bounds on $\sigma(G)$, and conjectured that $\sigma(T)=1$ for every tree $T$ on more than one vertex. Bounds on $\sigma(G)$ in terms of an optimal clique cover for $G$ were given by Gould and R\"{o}dl \cite{GR91}, and Harary's conjecture for trees was proved by Ellingham \cite{Ell93}. There has since been a great deal of work on the sum number and variations thereof; see \cite[Section 7.1]{Survey} for (much) more detail, including on the exclusive sum number defined below.

The exclusive sum number was introduced by Nagamochi, Miller and Slamin \cite{NMS01}, and further developed by Miller, Patel, Ryan, Sugeng, Slamin and Tuga \cite{MPRSST}. Given a sum labelling $f$ of a graph, a vertex $w$ is a \textit{working vertex} if $f(w)=f(u)+f(v)$ for some edge $uv$. If $G$ is a graph with sum number $r$ and $f$ is a sum labelling of $G\cup\overline{K}_r$, then necessarily all vertices not in $G$ are working vertices, since otherwise one of them could be omitted to obtain a sum labelling of $G\cup\overline{K}_{r-1}$. A sum labelling $f$ is an \textit{exclusive sum labelling} if these are the only working vertices, and the \textit{exclusive sum number}, $\epsilon(G)$, of $G$ to be the minimum $r$ for which $G\cup\overline{K}_r$ has an exclusive sum labelling. (We will generally assume $G$ to be a connected graph on more than one vertex, so that it will be clear which vertices may be working.) For further results on this parameter, see \cite{Ryan,MRR17}.

Before describing our result involving the exclusive sum number, we first give an equivalent, but more natural, definition. For two finite sets of positive integers $S,T$ we define $G_+(S,T)$ to be the graph on vertex set $S$, with two vertices being adjacent if and only if their sum is in $T$.
\begin{proposition}\label{prop:excl}For any finite graph $G$, $\epsilon(G)$ is the smallest possible value of $\abs{T}$ 
among pairs $S,T$ satisfying $G\cong G_+(S,T)$.
\end{proposition}
\begin{proof}
Note that if $f\colon V(G\cup \overline K_r)\to\zp$ is an exclusive sum labelling for $G$, then $G\cong G_+(S,T)$ where $S=f(V(G))$ and $T=f(V(\overline K_r))$ are disjoint. Conversely, if $G\cong G_+(S,T)$ for some pair of disjoint sets $S,T$ then we may recover an exclusive sum labelling $f\colon G\cup\overline K_{\abs{T}}\to\zp$ from this isomorphism: set $f(v)=\phi(v)$ if $v\in V(G)$, where $\phi\colon V(G)\to S$ is the isomorphism, and choose the remaining values of $f$ such that $f(V(\overline K_{\abs T}))=T$. 

Consequently $\epsilon(G)$ is the smallest value of $\abs{T}$ among pairs of disjoint sets $S,T$ such that $G\cong G_+(S,T)$. However, we may equivalently omit the disjointness condition, since $G_+(S,T)\cong G_+(S+\{k\},T+\{2k\})$ for any $k\in\zp$, and, for any fixed finite sets $S$ and $T$, choosing $k$ sufficiently large ensures $S+\{k\}$ and $T+\{2k\}$ are disjoint.\end{proof}

Harrington et al.\ \cite[Theorem 2.5]{HHKRW2} observed that the exclusive sum number gives an upper bound for the sum index, and proved that there is strict inequality for one specific graph. (There are also graphs for which equality occurs; for example $\sm[K_n]=\epsilon(K_n)$.) Motivated by this, they conjecture that the difference can be made arbitrarily large.
\begin{conjecture}\label{exclusive}For any positive integer $N$ there exists a graph $G$ with $\epsilon(G)-\sm>N$.\end{conjecture}
In fact Conjecture \ref{exclusive} is true in a strong sense, and this follows immediately from known results on the sum number. Clearly $\epsilon(G)\geq\sigma(G)$, and it was shown by Gould and R\"{o}dl \cite{GR91} that there exist graphs with $\sigma(G)=\Theta(\abs{G}^2)$, and by Nagamochi, Miller and Slamin \cite{NMS01} that such graphs are quite common. Since trivially $\sm=O(\abs{G})$, such graphs satisfy $\sm=o(\epsilon(G))$. However, both proofs referred to are non-constructive, and in fact no specific construction is known for which $\sigma(G)=\Omega(\abs{G})$ (see \cite{MPRSST,Survey}). In Section \ref{sec:exclusive} we give a constructive proof of Conjecture \ref{exclusive} by giving an explicit sequence of graphs with unbounded exclusive sum index for which the two parameters differ by a constant multiplicative factor.

We conclude by giving some natural open questions suggested by this work. In what follows, we often use the notation $[m]:=\{1,\ldots,m\}$ for $m\in\zp$.

\section{Graphs with large sum index and bounded difference index}\label{high-sum}
In this section we show that the sum index cannot be bounded by a function of the difference index, contradicting Conjecture \ref{conj:sum-dif}. If $G$ satisfies $\df<\ceil{\sm}$ then $G$ cannot be bipartite by \cite[Theorem 3.3]{HHKRW2}. However, we can give triangle-free examples, and indeed can forbid short odd cycles of arbitrary length.

We start by giving a new lower bound on the sum index.
\begin{theorem}\label{thm:long-cycles}If $G$ contains $s$ cycles of length $2k+1$ then $\sm\geq\sqrt[2k+1]{(4k+2)s}+1$.\end{theorem}
\begin{proof}Let $f:V(G)\to\z$ be an injective vertex labelling chosen such that $f_+(e)$ takes $\sm$ different values. Fix a cycle $v_1\cdots v_{2k+1}$ of length $2k+1$, and for ease of notation write $v_{2k+2}=v_1$. Note that \[2f(v_1)=\sum_{i=1}^{k+1}f_+(v_{2i-1}v_{2i})-\sum_{i=1}^{k}f_+(v_{2i}v_{2i+1}),\]
	and similar expressions for the other vertices may be obtained by symmetry.
	Thus, given the ordered edge labels around a $(2k+1)$-cycle we may deduce the ordered vertex labels. The ordered edge labels around a $(2k+1)$-cycle have the property that no two consecutive labels are equal and all labels are taken from the set $f_+(E(G))$ of $\sm$ elements. Thus every $(2k+1)$-tuple of ordered edge labels corresponds to a proper $\sm$-colouring of $C_{2k+1}$. Recall that the chromatic polynomial of $C_{2k+1}$ is given by $P(\lambda)=(\lambda-1)^{2k+1}-(\lambda-1)$, and so there are fewer than $(\sm-1)^{2k+1}$ possible tuples arising as the ordered edge labels of some $(2k+1)$-cycle. However, there are $s$ cycles of this length, each gives $4k+2$ tuples of ordered edge labels (since the starting point and direction may be chosen arbitrarily), and all these tuples must be different. Consequently
	\[(\sm-1)^{2k+1}\geq(4k+2)s,\]
	from which the claimed bound follows.
\end{proof}
As a consequence we deduce the following.
\begin{theorem}\label{thm:large-sum}There exist connected graphs with difference index $2$ but arbitrarily large sum index and girth.\end{theorem}
\begin{figure}
	\centering
	\begin{tikzpicture}
		\draw[fill] (0,0) circle (0.05) node[anchor=north]{$1$};
		\draw[fill] (2.5,0) circle (0.05) node[anchor=north]{$2k+1$};
		\draw[fill] (5,0) circle (0.05) node[anchor=north]{$4k+1$};
		\draw[fill] (7.5,0) circle (0.05) node[anchor=north]{$2(s-1)k+1$};
		\draw[fill] (10,0) circle (0.05) node[anchor=north]{$2sk+1$};
		\draw[fill] (0.4,1.5) circle (0.05) node[anchor=south]{$2$};
		\draw[fill] (0.8,1.5) circle (0.05) node[anchor=south]{$3$};
		\draw[fill] (2.1,1.5) circle (0.05) node[anchor=south]{$2k$};
		\draw[fill] (2.9,1.5) circle (0.05);
		\draw[fill] (3.3,1.5) circle (0.05);
		\draw[fill] (4.6,1.5) circle (0.05) node[anchor=south]{$4k$};
		\draw[fill] (7.9,1.5) circle (0.05);
		\draw[fill] (8.3,1.5) circle (0.05);
		\draw[fill] (9.6,1.5) circle (0.05) node[anchor=south]{$2sk$};
		\draw (0,0) -- (0.4,1.5) -- (0.8,1.5);
		\draw[dashed] (0.8,1.5) -- (2.1,1.5);
		\draw (0,0) -- (2.5,0) -- (2.1,1.5);
		\draw (2.5,0) -- (2.9,1.5) -- (3.3,1.5);
		\draw[dashed] (3.3,1.5) -- (4.6,1.5);
		\draw (2.5,0) -- (5,0) -- (4.6,1.5);
		\draw (7.5,0) -- (7.9,1.5) -- (8.3,1.5);
		\draw[dashed] (8.3,1.5) -- (9.6,1.5);
		\draw (7.5,0) -- (10,0) -- (9.6,1.5);
		\draw[dashed] (5,0) -- (7.5,0);
	\end{tikzpicture}
	\caption{Chaining $s$ odd cycles with difference index $2$.}\label{fig:long-cycles}
\end{figure}
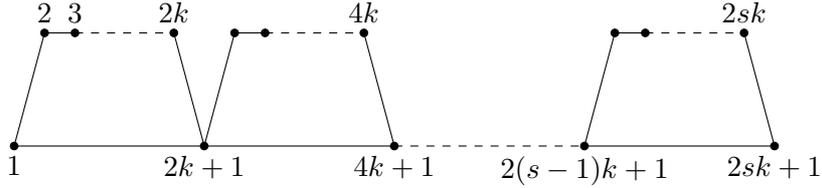
\begin{proof}Let $G$ consist of a path $v_1\cdots v_{2sk+1}$ together with the additional edges
\[\{v_{2(a-1)k+1}v_{2ak+1}\mid a\in[s]\};\]
this forms a chain of $s$ cycles, each of length $2k+1$. Let $f$ be the function $v_i\mapsto i$; see Figure \ref{fig:long-cycles} for $G$ with this labelling. Then $f_-:E(G)\to\{1,2k\}$, and since the only graphs with difference index $1$ are the path-forests we have $\df[T_n]=2$. However, $\sm[G]\geq \sqrt[2k+1]{(4k+2)s}+1$ and $G$ has girth $2k+1$; we may choose $s,k$ to make both parameters arbitrarily large.
\end{proof}

\section{Graphs with difference index larger than conjectured}\label{low-sum}
We now show that Conjecture \ref{conj:sum-dif} also fails in the other direction. In \cite[Theorem 3.1]{HHKRW2} it is observed that $\df\geq\delta(G)$. We give a strengthening of this fact, and establish a similar bound for the sum index. Write $\delta_k(G)$ for the $k$th smallest degree.
\begin{theorem}\label{delta2}For any graph $G$ we have $\df\geq\max_{k\geq 1}(\delta_{2k}(G)-k+1)$.
\end{theorem}
\begin{proof}
Let $G$ be a graph and choose $k$ to maximise the expression $\delta_{2k}(G)-k+1$. Fix an injective vertex labelling $f$, and order the vertices
$v_1,\ldots,v_n$ such that $f(v_i)<f(v_j)$ if and only if $i<j$. Observe
that at least one of $\{v_1,\ldots,v_k,v_{n-k+1},\ldots,v_n\}$ has degree at least $\delta_{2k}(G)$. Suppose without loss of generality that $d(v_i)\geq\delta_{2k}(G)$ for some $i\leq k$. Now $v_i$ has at least $d(v_i)-k+1$ neighbours $v_j$ with $j>i$, and hence $f(v_j)>f(v_i)$; the edge to such a neighbour gets label $f(v_j)-f(v_i)$, giving at least $\delta_{2k}(G)-k+1$ distinct labels.
\end{proof}
\begin{theorem}\label{delta-sum}For any graph $G$ we have $\sm\geq\max_{k\geq 1}(\delta_{k}(G)+\delta_{k+1}(G)-k)$.
\end{theorem}
\begin{remark}This gives a stronger lower bound for regular graphs than any previously known, since by taking $k=1$ we obtain $\sm\geq 2d-1$ if $G$ is $d$-regular. In \cite{HHKRW2} it is observed that $\sm\geq\chi'(G)$. However, by Vizing's theorem, $\chi'(G)\leq d+1$ if $G$ is $d$-regular, and so there is a considerable gap between these two bounds for large $d$.\end{remark}
\begin{proof}
Let $G$ be a graph and choose $k$ to maximise the expression $\delta_{k}(G)+\delta_{k+1}(G)-k$. Fix an injective vertex labelling $f$. Let $V_k$ be the $k-1$ vertices of smallest degree, and let $u,v\in V(G)\setminus V_k$ be chosen so that $f(u)$ is as small as possible and $f(v)$ is as large as possible.

Write $V_k^-$ for the set $\{w\in V_k:f(w)<f(u)\}$ and $V_k^+$ for the set $\{w\in V_k:f(w)>f(v)\}$. Note that every edge of the form $uw$ where $w\not\in V_k^+$ has sum label at most $f(u)+f(v)$, and all these labels are distinct, and likewise every edge of the form $vw$ where $w\not\in V_k^-$ has sum label at least $f(u)+f(v)$, and all these labels are distinct. The only label which may appear in both lists is $f(u)+f(v)$, if $uv$ is an edge. Consequently, combining these gives a set of at least $d(u)-\abs{V_k^+}+d(v)-\abs{V_k^-}-1$ edges with distinct labels. This is at least $\delta_{k}(G)+\delta_{k+1}(G)-\abs{V_k}-1$, as required.
\end{proof}
We immediately obtain a new result: the precise sum index of the prism graphs. The difference index of these graphs was already known \cite[Theorem 3.1, Corollary 3.13]{HHKRW2}.
\begin{corollary}For each $n\geq 3$, the prism graph $\Pi_n=C_n\square K_2$ satisfies $\sm[\Pi_n]=5$.
\end{corollary}
\begin{proof}Since $\Pi_n$ is $3$-regular, Theorem \ref{delta-sum} gives $\sm[\Pi_n]\geq 5$, but also $\sm[\Pi_n]\leq 5$ by \cite[Corollary 2.13]{HHKRW2}.
\end{proof}

We can now give a simple example where the difference index exceeds the conjectured bound. Let $\hat{K}_n$ be the complete graph on $n\geq 4$ vertices with one edge subdivided.
\begin{corollary}For each $n\geq 4$ we have $\df[\hat{K}_n]=\ceil{\sm[\hat{K}_n]/2}+1$.\end{corollary}
\begin{proof}Let $w$ be the new vertex of degree $2$, and let $u$ and $v$ be the two neighbours of $w$.

By Theorem \ref{delta2}, with $k=1$, we have $\df[\hat{K}_n]\geq n-1$, and this is an equality as shown by any labelling $f_1:V(\hat{K}_n)\to\{0,\ldots,n\}$ with $f_1(u)=0$ and $f_1(v)=n$, for which the difference labelling takes values $1,\ldots,n-1$.

By Theorem \ref{delta-sum}, with $k=2$, we have $\sm[\hat{K}_n]\geq 2n-4$, and this is also an equality as shown by any labelling $f_2:V(\hat{K}_n)\to\{0,\ldots,n\}$ with $f_2(u)=n-1$, $f_2(v)=n$ and $f_2(w)=0$, for which the sum labelling takes values $3,\ldots, 2n-2$. Consequently $\sm[\hat{K}_n]/2=n-2=\df[\hat{K}_n]-1$.\end{proof}

We will improve this construction to give graphs where the difference index exceeds the conjectured bound by an arbitrarily large amount, although the relative error will be small. Let $\hat{K}_n^{(k,k)}$ be the graph obtained from the complete graph on $n\geq 2k$ vertices $v_1,\ldots,v_n$ by subdividing all edges $v_iv_j$ where $i,j\leq k$ and all edges $v_iv_j$ where $i,j>n-k$ (that is, we subdivide the edges of two disjoint cliques of order $k$). Figure \ref{fig:pentagon} shows the case $k=2$, $n=2k+\binom k2=5$. The motivation for this construction is that a vertex labelling $f$ for $K_n$ achieves the sum index if and only if its values are in arithmetic progression, and the same is true for the difference index. 
For such a labelling, the $2k$ least frequent labels of $f^+$ occur on the edges of two disjoint cliques 
(whose vertices have the $k$ highest and $k$ lowest labels), whereas the $2k$ least frequent labels of 
$f^-$ occur on a subset of the edges of a $K_{k,k}$ (with the same vertices). 
In what follows, we refer to the vertices $v_1,\ldots,v_k$ as the \textit{branch vertices}.
\begin{figure}
\centering
\begin{tikzpicture}
\draw[fill] (0,0) circle (0.05) node[anchor=north]{$1$};
\draw[fill] (2,0) circle (0.05) node[anchor=north]{$5$};
\draw[fill] (2.618,1.902) circle (0.05) node[anchor=west]{$4$};
\draw[fill] (-.618,1.902) circle (0.05) node[anchor=east]{$2$};
\draw[fill] (1,3.078) circle (0.05) node[anchor=south]{$3$};
\draw[fill] (-.618,.851) circle (0.05) node[anchor=east]{$6$};
\draw[fill] (2.618,.851) circle (0.05) node[anchor=west]{$0$};
\draw (0,0) -- (1,3.078) -- (2,0) -- (-.618,1.902) -- (2.618,1.902) -- cycle;
\draw (0,0) -- (-.618,.851) -- (-.618,1.902) -- (1,3.078) -- (2.618,1.902) -- (2.618,.851) -- (2,0) -- cycle;
\end{tikzpicture}
\caption{The graph $\hat{K}_5^{(2,2)}$, with a labelling demonstrating that its sum index is $5$.}\label{fig:pentagon}
\end{figure}
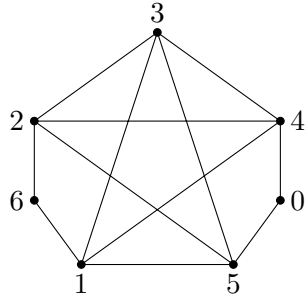
\begin{theorem}\label{knkk}For each $k\geq 2$ and $n\geq\binom{k}{2}+2k$, we have 
\[\smb= 2n-2k-1\]
and
\[\dfb\geq n-2k/3-1=\bceil{\smb/2}+k/3-1.\]\end{theorem}
\begin{proof}For the sum index, we first prove the lower bound. Let $f_1$ be any vertex labelling, and let $i$ and $j$ be the indices with $f_1(v_i)$ as small as possible and $f_1(v_j)$ as large as possible. For the corresponding sum labelling, there are at least $n-k$ edges from $v_i$ to other branch vertices, each of which have different labels of at most $f_1(v_i)+f_1(v_j)$, and there are at least $n-k-1$ edges from $v_j$ to other branch vertices which are not $v_i$, each of which have different labels exceeding $f_1(v_i)+f_1(v_j)$. Thus $\smb\geq 2n-2k-1$. If $f_1$ is chosen such that $f_1(v_i)=i$, the subdivision vertices between two vertices of small index receive labels $n+1,\ldots,\binom k2$, and the subdivision vertices between two vertices of high index receive labels $1-\binom k2,\ldots,0$ then, since $n\geq \binom k2+2k$, the sum labelling uses only labels from $\{k+2,\ldots,2n-k\}$, giving the corresponding upper bound. See Figure \ref{fig:pentagon} for an example labelling.

Now suppose $\dfb< n-2k/3-1$, and let $f_2$ be a vertex labelling achieving the difference index. Write $S=\{k+1,\ldots,n-k\}$, and $T=[n]\setminus S$, i.e.\ a branch vertex $v_i$ is adjacent to all other branch vertices if and only if $i\in S$. The condition on $n$ ensures that $S$ is nonempty. Choose $i\in S$ such that $f_2(v_i)$ is as small as possible, and $j\in[n]$ so that $f_2(v_j)$ is as large as possible. Set 
\begin{align*}A&=\{a\in[n]:f_2(v_a)<f_2(v_i)\},\\B&=\{b\in[n]:f_2(v_b)>f_2(v_i)\}.\end{align*}
Since $i\in S$ there is an edge $v_iv_b$ for each $b\in B$, and all these edges get different labels which are at most $f_2(v_j)-f_2(v_i)$. For each $a\in A$ such that $v_jv_a$ is an edge, $v_jv_a$ gets a different label which is greater than $f(v_j)-f(v_i)$. Since $\abs{A}+\abs{B}=n-1$, there must be more than $2k/3$ values $a\in A$ for which $v_jv_a$ is not an edge; i.e., without loss of generality, $j\in[k]$ and $\babs{A\cap[k]}>2k/3$.

Now choose $h\in [n]\setminus[k]$ so that $f_2(v_h)$ is as large as possible. Set 
\[C=\{b\in B:f_2(v_b)\leq f_2(v_h)\}.\]
This choice of $h$ ensures that $B\setminus C\subseteq B\cap[k]$, and since $\babs{A\cap[k]}>2k/3$ we have $\babs{B\cap[k]}<k/3$, meaning that $\abs{A}+\abs{C}> n-1-k/3$. Now for each $c\in C$ the edge $v_iv_c$ has label at most $f_2(v_h)-f_2(v_i)$, and for each $a\in A\cap[n-k]$ there is an edge $v_hv_a$ with label exceeding $f_2(v_h)-f_2(v_i)$, with all these labels being distinct. Since fewer than $n-2k/3-1$ labels are used by assumption, we have 
\begin{align*}n-2k/3-1&>\babs{A\cap[n-k]}+\abs{C}\\
&=\abs{A}+\abs{C}-\babs{A\cap\{n-k+1,\ldots, n\}}\\
&> n-1-k/3-\babs{A\cap\{n-k+1,\ldots, n\}},\end{align*}
and so $\babs{A\cap \{n-k+1,\ldots,n\}}>k/3$. Consequently $\abs{A\cap T}>k$.

Now choose $i'\in S$ so that $f_2(v_{i'})$ is as large as possible, and set
\[A'=\{a\in[1,n]:f_2(v_a)>f_2(v_{i'})\}.\]
By a similar argument, we have $\abs{A'\cap T}>k$. Since $A$ and $A'$ are disjoint, this gives 
\[\babs{(A\cup A')\cap T)}>2k=\abs T,\]
a contradiction.
\end{proof}
Choosing the smallest value of $n$ for each $k$ in Theorem \ref{knkk}, we have $\smb=\Theta(k^2)$ and so $\dfb=\bceil{\smb/2}+\Omega\bigl(\smb^{1/2}\bigr)$.
\section{Graphs with sum index far from exclusive sum number}\label{sec:exclusive}
In this section we give a constructive proof of Conjecture \ref{exclusive} by giving explicit graphs $G_N$ satisfying $\epsilon(G_N)-\sm[G_N]\geq N$ for any fixed $N$; in fact $\abs{G_N}=O(N)$ so the two invariants differ by a constant proportion. Recall, however, that a much stronger form of Conjecture \ref{exclusive} follows from known, but non-constructive, results.

Our proof will use a stability result from additive combinatorics. For sets of integers (all sets we consider will be finite and nonempty) $A$ and $B$ we write $A+B$ for the \textit{sumset} $\{a+b\mid a\in A, b\in B\}$. It is easy to verify that $\abs{A+B}\geq\abs{A}+\abs{B}-1$, with equality if and only if $A$ and $B$ are arithmetic progressions with the same common difference. A natural question is whether this structural constraint is robust when replacing equality with near-equality. Stanchescu \cite{Sta96}, building on earlier work by Fre\u{\i}man \cite{Fre62} and by Lev and Smeliansky \cite{LS95}, showed that it is robust in the following sense. For a nonempty finite set $X$ of integers, write $\ell(X)=\max(X)-\min(X)$. Given two such sets $A$ and $B$, let $d(A,B)$ be the greatest common divisor of $\{a-a'\mid a,a'\in A\}\cup\{b-b'\mid b,b'\in B\}$. We write $\delta_{ij}$ to mean the Kronecker delta.
\begin{theorem}[{\cite[Corollary 2]{Sta96}}]\label{stability}Let $A$ and $B$ be two finite sets of integers.
If $\abs{A+B} < \abs A+\abs B+ \min(\abs A,\abs B)-2-\delta_{\ell(A)\ell(B)}$ then $\ell(A)/d(A,B) \leq \abs{A+B}-\abs{B}$ and $\ell(B)/d(A,B) \leq\abs{A+B}-\abs{A}$.
\end{theorem}
Note that, writing $d=d(A,B)$, we may find arithmetic progressions $A'$ and $B'$ of order $1+\ell(A)/d$ and $1+\ell(B)/d$ respectively, each having difference $d$, such that $A$ is contained in $A'$ and $B$ is contained in $B'$. Theorem \ref{stability} implies that if $A+B$ is close to minimal size then $A'$ and $B'$ are also not much larger than $A$ and $B$.

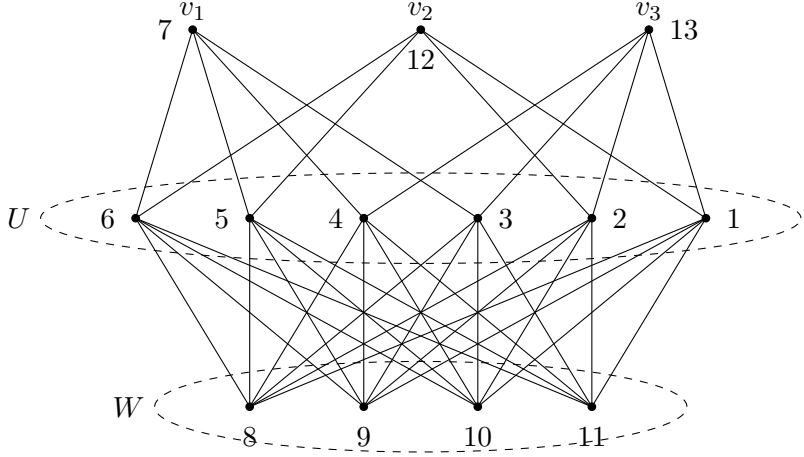
\begin{figure}
\centering
\begin{tikzpicture}
\draw[fill] (.75,2.5) circle (0.05) node[anchor=south]{$v_1$} node[label=left:$7$]{};
\draw[fill] (3.75,2.5) circle (0.05) node[anchor=south]{$v_2$} node[label=below:$12$]{};
\draw[fill] (6.75,2.5) circle (0.05) node[anchor=south]{$v_3$} node[label=right:$13$]{};
\foreach \x in {1,2,3}{\draw[fill] (9-1.5*\x,0) circle (0.05) node[label=right:$\x$]{};}
\foreach \x in {4,5,6}{\draw[fill] (9-1.5*\x,0) circle (0.05) node[label=left:$\x$]{};}
\foreach \x in {8,...,11}{\draw[fill] (1.5*\x-10.5,-2.5) circle (0.05) node[label=below:$\x$]{};}
\draw[dashed] (3.75,0) circle [x radius=5, y radius=0.6];
\node at (-1.25,0) [anchor=east] {$U$};
\draw[dashed] (3.75,-2.5) circle [x radius=3.5, y radius=0.6];
\node at (.25,-2.5) [anchor=east] {$W$};
\foreach \x in {0,1.5,...,7.5}{%
\foreach \y in {1.5,3,4.5,6}{\draw (\x,0) -- (\y,-2.5);}}
\foreach \x in {0,1.5,3,4.5}{\draw (\x,0) -- (.75,2.5);}
\foreach \x in {0,1.5,6,7.5}{\draw (\x,0) -- (3.75,2.5);}
\foreach \x in {3,4.5,6,7.5}{\draw (\x,0) -- (6.75,2.5);}
\end{tikzpicture}
\caption{The graph $G_{6,4}$ with a labelling demonstrating $\sm[G_{6,4}]\leq 10$.}\label{fig:a+b}
\end{figure}

Next we give the construction used in our proof. 
Fix positive integers $k\geq 3$ and $n\geq 3k-6$. Take a complete bipartite graph $K_{n,k}$ and let $U$ and $W$ be the classes of order $n$ and $k$ respectively. Add three additional vertices $v_1,v_2,v_3$, each adjacent to exactly $n-k+2$ vertices of $U$, in such a way that the three sets given by $U\setminus N(v_i)$ are disjoint. This is possible since each set has size $k-2$ and $n\geq 3k-6$. Denote the resulting graph by $G_{n,k}$. Figure \ref{fig:a+b} shows the case $k=4$ and $n=3k-6=6$.
\begin{theorem}The graph defined above satisfies $\sm[G_{n,k}]\leq n+k$ and $\epsilon(G_{n,k})\geq n+2k-3$.
\end{theorem}
\begin{remark}In particular, taking $n=3k-6$ gives a sequence of graphs with unbounded sum index satisfying $\sm\leq\frac{4\epsilon(G)+6}{5}$.
\end{remark}
\begin{proof}The first statement arises from a bijective labelling $f:V(G_{k,n})\to[n+k+3]$ satisfying the following conditions.
\begin{itemize}
\item $f(v_1)=n+1$ and $f(v_3)=n+k+3$;
\item $f(u)=1$ for some $u\in U\setminus N(v_1)$;
\item $f(u')=n$ for some $u'\in U\setminus N(v_3)$;
\item and $f(U)=[n]$.
\end{itemize}
In this case the smallest value of $f_+(e)$ for $e\in E(G_{k,n})$ will be $n+2$ and the largest will be $2n+k+2$. See Figure \ref{fig:a+b} for an example.

Now we bound the exclusive sum number using Proposition \ref{prop:excl}. Suppose we can express $G_{n,k}$ as $G_+(S,T)$ for some sets $S,T$ with $\abs T\leq n+2k-4$. Fix an isomorphism $\phi:V(G_{n,k})\to S$, write $A=\phi(U)$ and $B=\phi(W)$, and set $d=d(A,B)$. By definition of $G_{n,k}$ and $G_+(S,T)$ we must have $a+b\in T$ for each $a\in A, b\in B$. 

We have $\abs{A+B}\leq\abs{T}=n+2k-4$; write $t=n+2k-4-\abs{A+B}$. Since $\abs{A}=n$, $\abs{B}=k$, and the conditions on $k$ and $n$ imply $k\leq n$, we have $\abs{A+B}<\abs{A}+\abs{B}+\min\{\abs{A},\abs{B}\}-3$. Consequently, by Theorem \ref{stability} we have $\ell(A)/d\leq n+k-4-t$ and $\ell(B)/d\leq 2k-4-t$. Also, since $\ell(A+B)=\ell(A)+\ell(B)$, we have \[\ell(A+B)/d\leq 2\abs{A+B}-n-k\leq\abs{A+B}+k-4-t.\]
It follows that there are arithmetic progressions $A'\supseteq A$, $B'\supseteq B$ and $C\supseteq A+B$, all with difference $d$, such that each of $A'\setminus A$, $B'\setminus B$ and $C\setminus(A+B)$ has size at most $k-3-t$. All elements of $B'$ are congruent modulo $d$; write $B^*$ for the set of all integers congruent to elements of $B'$ modulo $d$.

Write $r_i=\phi(v_i)$ for each $i\in [3]$. We claim that each $r_i\in B^*\setminus B'$. First, if $r_i\not\in B^*$ then we have $r_i+a\not\in A+B$ for each $a\in A$. For $n-(k-2)$ choices of $a$, the vertex associated with $a$ is a neighbour of $v_i$, meaning $r_i+a\in T$. This would imply $\abs{T}\geq\abs{A+B}+n-k+2\geq 2n+1$, a contradiction, so we must have $r_i\in B^*$. Secondly, if $r_i\in B'$ then $r_i+a\in C$ for each $a\in A$, and since $\abs{C\setminus (A+B)}\leq k-3-t$ there are at most $k-3-t$ choices of $a\in A$ for which $r_i+a\not\in T$. But there are $k-2$ vertices in $U$ not adjacent to $v_i$, again giving a contradiction. So $r_i\in B^*\setminus B'$.

Now, for each $i$ either $r_i>\max(B')$ or $r_i<\min(B')$, and so one of these inequalities must hold for at least two of the three possible values of $i$; assume without loss of generality that $r_1,r_2>\max(B')$. Note that there must be exactly $k-2$ choices of $a\in A$ such that $a+r_1\not\in T$. There are at most $k-3-t$ choices of $a\in A$ for which $a+r_1\in C\setminus(A+B)$, so there are at least $t+1$ choices with $a+r_1\not\in C\cup T$. Also, since $T\setminus C\subseteq T\setminus (A+B)$, there are at most $t$ choices $a\in A$ such that $a+r_1\in T\setminus C$. Note that since $r_1>\max(B')$, the $a\in A$ with $a+r_1\not\in C$ are precisely the $s$ largest elements of $A$ for some $s$. Since these $s$ contain at least $t+1$ choices of $a$ with $a+r_1\not\in T$ and at most $t$ choices of $a$ with $a+r_1\in T$, it follows that at least $t+1$ of the $2t+1$ largest elements of $A$ satisfy $a+r_1\not\in T$, i.e.\ are images under $\phi$ of vertices of $U$ which are not adjacent to $v_1$. However, similarly, at least $t+1$ of these $2t+1$ elements are images under $\phi$ of vertices of $U$ which are not neighbours of $v_2$, and so at least one must be $\phi(u)$ for some $u\in U$ not adjacent to $v_1$ or $v_2$, giving a contradiction.
\end{proof}
\section{Open questions}
These results raise several natural questions. Notice that all the counterexamples we obtain to Conjecture \ref{conj:sum-dif} are non-bipartite. For Section \ref{high-sum} this is necessarily the case, since $\df\geq\ceil{\sm/2}$ for bipartite graphs \cite{HHKRW}. In the original version of this paper we asked whether a bipartite counterexample to Conjecture \ref{conj:sum-dif} in the other direction exists. In \cite{HHKRW2} (the updated version of \cite{HHKRW}), a counterexample is given which not only is bipartite but in fact is a tree. However, this example only has a discrepancy of $1$.
\begin{question}If $G$ is bipartite, can $\df-\ceil{\sm/2}$ be arbitrarily large?\end{question}
We would also be interested to know by how far $\df$ can exceed $\ceil{\sm/2}$ in general. In Section \ref{low-sum} we showed that the absolute discrepancy can be arbitrarily large, but the relative discrepancy of our examples is small. Is it possible to improve our results to give a constant relative discrepancy?
\begin{question}\label{q2}Is there a constant $\alpha>1$ for which there are graphs with arbitrarily large sum index satisfying $\df\geq(1+\alpha)\ceil{\sm/2}$?\end{question}
We conjecture that, even if the answer to Question \ref{q2} were positive, the results of Section \ref{low-sum} could not be improved any further than that. Recall that Harrington et al.\ proved that no bipartite graph has larger difference index than sum index \cite{HHKRW2}; we expect this relationship to hold more generally.
\begin{conjecture}Every graph $G$ satisfies $\df\leq\sm$.\end{conjecture}
Finally, we turn to the results of Section \ref{sec:exclusive}. By non-constructive arguments mentioned in Section \ref{intro}, we know there are graphs with $\sm=o(\epsilon(G))$. However, we cannot currently rule out the possibility that much more is true in the sense of Theorem \ref{thm:large-sum}.
\begin{question}Do there exist graphs with arbitrarily large exclusive sum number and bounded sum index?\end{question}

\end{document}